\documentclass[10pt]{article}         
\usepackage{authblk,graphicx,amsmath,amsfonts,amssymb,amscd,amsmath,breqn}  
\usepackage{amsthm}

\DeclareMathOperator{\csch}{csch}

\title{Table in Gradshteyn and Ryzhik: Derivation of definite integrals of a Hyperbolic Function}

\author[1]{Robert Reynolds}
\author[2]{Allan Stauffer}
\affil[1]{Faculty of Science, York University}
\affil[2]{Department of Mathematics, York University}
\affil[ ]{milver@my.yorku.ca, stauffer@yorku.ca}

\begin{document}

\maketitle

\begin{abstract}
We present a method using contour integration to derive definite integrals and their associated infinite sums which can be expressed as a special function. We give a proof of the basic equation and some examples of the method. The advantage of using special functions is their analytic continuation which widens the range of the parameters of the definite integral over which the formula is valid. We give as examples definite integrals of logarithmic functions times a trigonometric function. In various cases these generalizations evaluate to known mathematical constants such as Catalan's constant and $\pi$.

\end{abstract}

\section{Introduction}

We will derive integrals as indicated in the abstract in terms of special functions. Some special cases of these integrals have been reported in Gradshteyn and Ryzhik  \cite{grad}. In 1867 David Bierens de Haan \cite{bdh} derived hyperbolic integrals of the form 
 \begin{equation*}   
 \int_{0}^{\infty}\frac{\sinh (a x) \left(e^{-m x} (\log (\alpha )-x)^k-e^{m x} (\log (\alpha )+x)^k\right)}{(\cosh (a x)+\cos (t))^2}dx
 \end{equation*}

 In our case the constants in the formulas are general complex numbers subject to the restrictions given below.  The derivations follow the method used by us in \cite{reyn3}. The generalized Cauchy's integral formula is given by

\begin{equation}\label{intro:cauchy}
\frac{x^k}{k!}=\frac{1}{2\pi i}\int_{C}\frac{e^{wx}}{w^{k+1}}dw.
\end{equation}

 This method involves using a form of equation (\ref{intro:cauchy}) then multiplys both sides by a function, then takes a definite integral of both sides. This yields a definite integral in terms of a contour integral. Then we multiply both sides of equation (\ref{intro:cauchy})  by another function and take the infinite sum of both sides such that the contour integral of both equations are the same.

\section{Derivation of the definite integral of the contour integral}

We use the method in \cite{reyn3}. Here the contour is similar to Figure 2 in \cite{reyn3}. Using a generalization of Cauchy's integral formula we first replace $x$ by $ix+\log(a)$ then multiply both sides by $e^{mx}$ for the first equation and the replace $x$ with $-x$ and multiplying both sides by $e^{-mx}$ to get the second equation. Then we subtract these two equations, followed by multiplying both sides by $-\frac{\sinh (a x)}{2 (\cosh (a x)+\cos (t))^2}$ to get

\begin{multline}\label{albania:eq1}
-\frac{\sinh (a x) \left(e^{-m x} (\log (\alpha )-x)^k-e^{m x} (\log (\alpha )+x)^k\right)}{2 k! (\cosh (a x)+\cos (t))^2}\\\\
=\frac{1}{2\pi i}\int_{C}\frac{w^{-k-1} \alpha ^w \sinh (a x) \sinh (x (m+w))}{(\cosh (a x)+\cos (t))^2}dw
\end{multline}

where the logarithmic function is defined in equation (4.1.2) in \cite{as}. We then take the definite integral over $x \in [0,\infty)$ of both sides to get

\begin{dmath}\label{derive:eq1}
-\int_{0}^{\infty}\frac{\sinh (a x) \left(e^{-m x} (\log (\alpha )-x)^k-e^{m x} (\log (\alpha )+x)^k\right)}{2 k! (\cosh (a x)+\cos (t))^2}dx
=\frac{1}{2\pi i}\int_{0}^{\infty}\int_{C}\frac{w^{-k-1} \alpha ^w \sinh (a x) \sinh (x (m+w))}{(\cosh (a x)+\cos (t))^2}dwdx
=\frac{1}{2\pi i}\int_{C}\int_{0}^{\infty}\frac{w^{-k-1} \alpha ^w \sinh (a x) \sinh (x (m+w))}{(\cosh (a x)+\cos (t))^2}dxdw
=\frac{1}{2\pi i}\int_{C}\frac{\pi  m w^{-k-1} \csc (t) \alpha ^w \csc \left(\frac{\pi  (m+w)}{a}\right) \sin \left(\frac{t (m+w)}{a}\right)}{a^2}dw+\frac{1}{2\pi i}\int_{C}\frac{\pi  w^{-k} \csc (t) \alpha ^w \csc \left(\frac{\pi (m+w)}{a}\right) \sin \left(\frac{t (m+w)}{a}\right)}{a^2}dw
\end{dmath}

from equation (2.5.48.18) in \cite{prud} and the integrals are valid for $a$, $m$, $k$, $t$ and $\alpha$ complex and $-1<Re(w+m)<0$ and $Re(\alpha)\neq 0$. We are able to switch the order of integration over $w$ and $x$ using Fubini's theorem since the integrand is of bounded measure over the space $\mathbb{C} \times [0,\infty)$.

\section{Derivation of the infinite sum of the contour integral}

\subsection{Derivation of the first contour integral}

In this section we will again use the generalized Cauchy's integral formula to derive equivalent contour integrals. First we multiply equation (\ref{intro:cauchy}) by $e^{imt/\alpha}/2i$ then replace by $x$ by $p+it/\alpha$ for the first equation and then $p-it/\alpha$ for the second equation to get

\begin{equation}
\frac{i e^{-\frac{i m t}{a}} \left(\left(p-\frac{i t}{a}\right)^k-e^{\frac{2 i m t}{a}} \left(p+\frac{i t}{a}\right)^k\right)}{2 k!}=\frac{1}{2\pi}\int_{C}w^{-k-1} e^{w p} \sin \left(\frac{t
   (m+w)}{a}\right)dw
\end{equation}

Then we replace $p$ with $\pi i(2p+1)/a+\log(\alpha)$ and multiply both sides by $-\frac{2\pi }{a^1}$ to get

\begin{multline}\label{sum:eq1}
\frac{i e^{-\frac{i m t}{a}} \left(\left(\frac{i \pi  (2 p+1)}{a}-\frac{i t}{a}+\log (\alpha )\right)^k-e^{\frac{2 i m t}{a}} \left(\frac{i \pi  (2 p+1)}{a}+\frac{i t}{a}+\log (\alpha)\right)^k\right)}{2 k!}\\
 =\frac{1}{2\pi i}\int_{C}w^{-k-1} \sin \left(\frac{t (m+w)}{a}\right) e^{w \left(\log (\alpha )+\frac{i \pi  (2 p+1)}{a}\right)}dw
\end{multline}

Then we multiply both sides by $-\frac{2 i \pi}{a^2}   e^{\frac{i \pi  m (2 y+1)}{a}}$ and take the sum over $p\in[0,\infty)$ and simplify the left-hand side in terms of the Lerch function to get

\begin{multline}\label{derive:eq2}
\frac{2^k \pi ^{k+1} \left(\frac{i}{a}\right)^k e^{\frac{i m (\pi -t)}{a}} \left(\Phi \left(e^{\frac{2 i m \pi }{a}},-k,\frac{-t-i a \log (\alpha )+\pi }{2 \pi }\right)-e^{\frac{2 i m
   t}{a}} \Phi \left(e^{\frac{2 i m \pi }{a}},-k,\frac{t-i a \log (\alpha )+\pi }{2 \pi }\right)\right)}{a^2 k!}\\ 
   =\frac{1}{2\pi i}\sum_{p=0}^{\infty}\int_{C}w^{-k-1} \sin \left(\frac{t (m+w)}{a}\right) e^{w \left(\log (\alpha )+\frac{i \pi  (2 p+1)}{a}\right)}dw\\
   =\frac{1}{2\pi i}\int_{C}\sum_{p=0}^{\infty}w^{-k-1} \sin \left(\frac{t (m+w)}{a}\right) e^{w \left(\log (\alpha )+\frac{i \pi  (2 p+1)}{a}\right)}dw\\
   =\frac{1}{2\pi }\int_{C}\frac{\pi  w^{-k-1} \alpha ^w \csc \left(\frac{\pi  (m+w)}{a}\right) \sin
   \left(\frac{t (m+w)}{a}\right)}{a^2}dw
\end{multline}

from equation (1.232.3) in \cite{grad} where $\csch(ix)=-i\csc(x)$ from equation (4.5.10) in \cite{as} and $Im(w)>0$ for the sum to converge. The $\log$ terms cannot be combined in general.

\subsection{Derivation of the second contour integral}

Next we will derive the second equation by using equation (\ref{derive:eq2}), multiplying by $m\csc(t)$ and taking the infinite sum over $p\in[0,\infty)$ to get

\begin{multline}
\frac{2^k \pi ^{k+1} m \left(\frac{i}{a}\right)^k \csc (t) e^{\frac{i m (\pi -t)}{a}} \left(\Phi \left(e^{\frac{2 i m \pi }{a}},-k,\frac{-t-i a \log (\alpha )+\pi }{2 \pi
   }\right)-e^{\frac{2 i m t}{a}} \Phi \left(e^{\frac{2 i m \pi }{a}},-k,\frac{t-i a \log (\alpha )+\pi }{2 \pi }\right)\right)}{a^2 k!}\\
   =\frac{1}{2\pi i}\int_{C}\frac{\pi  m w^{-k-1} \csc (t) \alpha ^w \csc\left(\frac{\pi  (m+w)}{a}\right) \sin \left(\frac{t (m+w)}{a}\right)}{a^2}dw
\end{multline}

Then we replace $k$ with $k-1$ to get

\begin{multline}\label{derive:eq3}
\frac{2^{k-1} \pi ^k \left(\frac{i}{a}\right)^{k-1} \csc (t) e^{\frac{i m (\pi -t)}{a}} \left(\Phi \left(e^{\frac{2 i m \pi }{a}},1-k,\frac{-t-i a \log (\alpha )+\pi }{2 \pi
   }\right)-e^{\frac{2 i m t}{a}} \Phi \left(e^{\frac{2 i m \pi }{a}},1-k,\frac{t-i a \log (\alpha )+\pi }{2 \pi }\right)\right)}{a^2 (k-1)!}\\
   =\frac{1}{2\pi i}\int_{C}\frac{\pi  w^{-k} \csc (t) \alpha ^w \csc\left(\frac{\pi  (m+w)}{a}\right) \sin \left(\frac{t (m+w)}{a}\right)}{a^2}dw
\end{multline}

from equation (1.232.3) in \cite{grad} where $\csch(ix)=-i\csc(x)$ from equation (4.5.10) in \cite{as} and $Im(w)>0$ for the sum to converge.

\section{The Lerch function}

The Lerch function has a series representation given by
\begin{equation}\label{armenia:eq7}
\Phi(z,s,v)=\sum_{n=0}^{\infty}(v+n)^{-s}z^{n}
\end{equation}

where $|z|<1, v \neq 0,-1,..$ and is continued analytically by its integral representation given by

\begin{equation}\label{armenia:eq8}
\Phi(z,s,v)=\frac{1}{\Gamma(s)}\int_{0}^{\infty}\frac{t^{s-1}e^{-vt}}{1-ze^{-t}}dt=\frac{1}{\Gamma(s)}\int_{0}^{\infty}\frac{t^{s-1}e^{-(v-1)t}}{e^{t}-z}dt
\end{equation}

where $Re(v)>0$, or $|z| \leq 1, z \neq 1, Re(s)>0$, or $z=1, Re(s)>1$.  

\section{Definite integral in terms of the Lerch function} 

Since the right-hand sides of equation (\ref{derive:eq1}), (\ref{derive:eq2}) and (\ref{derive:eq3}) are equivalent we can equate the left-hand sides simplify the factorial to get

\begin{dmath}\label{definite:eq1}
\int_{0}^{\infty}\frac{\sinh (a x) \left(e^{-m x} (\log (\alpha )-x)^k-e^{m x} (\log (\alpha )+x)^k\right)}{(\cosh (a x)+\cos (t))^2}dx
=-\frac{k (2 \pi )^k
   \left(\frac{i}{a}\right)^{k-1} \csc (t) e^{\frac{i m (\pi -t)}{a}} \Phi \left(e^{\frac{2 i m \pi }{a}},1-k,\frac{-t-i a \log (\alpha )+\pi }{2 \pi }\right)}{a^2}+\frac{k (2 \pi )^k
   \left(\frac{i}{a}\right)^{k-1} \csc (t) e^{\frac{i m (\pi -t)}{a}+\frac{2 i m t}{a}} \Phi \left(e^{\frac{2 i m \pi }{a}},1-k,\frac{t-i a \log (\alpha )+\pi }{2 \pi }\right)}{a^2}-\frac{(2 \pi
   )^{k+1} m \left(\frac{i}{a}\right)^k \csc (t) e^{\frac{i m (\pi -t)}{a}} \Phi \left(e^{\frac{2 i m \pi }{a}},-k,\frac{-t-i a \log (\alpha )+\pi }{2 \pi }\right)}{a^2}+\frac{(2 \pi )^{k+1} m
   \left(\frac{i}{a}\right)^k \csc (t) e^{\frac{i m (\pi -t)}{a}+\frac{2 i m t}{a}} \Phi \left(e^{\frac{2 i m \pi }{a}},-k,\frac{t-i a \log (\alpha )+\pi }{2 \pi }\right)}{a^2}
\end{dmath}

The integral in equation (\ref{definite:eq1}) can be used as an alternative method to evaluating the Lerch function. 

\section{Evaluation of special cases of definite Integrals}

\subsection{Special case 1}

For this special case we will form a second equation using (\ref{definite:eq1}) by replacing $m$ by $-m$ taking the difference from the original equation and simplifying to get

\begin{dmath}\label{closedform:eq1}
-\int_{0}^{\infty}\frac{2 \sinh (a x) \sinh (m x) \left((\log (\alpha )-x)^k+(\log (\alpha )+x)^k\right)}{(\cosh (a x)+\cos (t))^2}dx=
\frac{k (2 \pi )^k \left(\frac{i}{a}\right)^{k-1} \csc (t) e^{-\frac{i m
   (\pi -t)}{a}} \Phi \left(e^{-\frac{2 i m \pi }{a}},1-k,\frac{-t-i a \log (\alpha )+\pi }{2 \pi }\right)}{a^2}-\frac{k (2 \pi )^k \left(\frac{i}{a}\right)^{k-1} \csc (t) e^{-\frac{i m (\pi
   -t)}{a}-\frac{2 i m t}{a}} \Phi \left(e^{-\frac{2 i m \pi }{a}},1-k,\frac{t-i a \log (\alpha )+\pi }{2 \pi }\right)}{a^2}-\frac{k (2 \pi )^k \left(\frac{i}{a}\right)^{k-1} \csc (t) e^{\frac{i
   m (\pi -t)}{a}} \Phi \left(e^{\frac{2 i m \pi }{a}},1-k,\frac{-t-i a \log (\alpha )+\pi }{2 \pi }\right)}{a^2}+\frac{k (2 \pi )^k \left(\frac{i}{a}\right)^{k-1} \csc (t) e^{\frac{i m (\pi
   -t)}{a}+\frac{2 i m t}{a}} \Phi \left(e^{\frac{2 i m \pi }{a}},1-k,\frac{t-i a \log (\alpha )+\pi }{2 \pi }\right)}{a^2}-\frac{(2 \pi )^{k+1} m \left(\frac{i}{a}\right)^k \csc (t) e^{-\frac{i
   m (\pi -t)}{a}} \Phi \left(e^{-\frac{2 i m \pi }{a}},-k,\frac{-t-i a \log (\alpha )+\pi }{2 \pi }\right)}{a^2}+\frac{(2 \pi )^{k+1} m \left(\frac{i}{a}\right)^k \csc (t) e^{-\frac{i m (\pi
   -t)}{a}-\frac{2 i m t}{a}} \Phi \left(e^{-\frac{2 i m \pi }{a}},-k,\frac{t-i a \log (\alpha )+\pi }{2 \pi }\right)}{a^2}-\frac{(2 \pi )^{k+1} m \left(\frac{i}{a}\right)^k \csc (t) e^{\frac{i m
   (\pi -t)}{a}} \Phi \left(e^{\frac{2 i m \pi }{a}},-k,\frac{-t-i a \log (\alpha )+\pi }{2 \pi }\right)}{a^2}+\frac{(2 \pi )^{k+1} m \left(\frac{i}{a}\right)^k \csc (t) e^{\frac{i m (\pi
   -t)}{a}+\frac{2 i m t}{a}} \Phi \left(e^{\frac{2 i m \pi }{a}},-k,\frac{t-i a \log (\alpha )+\pi }{2 \pi }\right)}{a^2}
\end{dmath}

\subsection{Special case 2}

For this special case we use equation (\ref{closedform:eq1}) setting $\alpha=1$ and taking the first partial derivative with respect to $m$ simplifying to get

\begin{dmath}\label{closedform:sp2}
\int_{0}^{\infty}\frac{x^k \sinh (a x) \cosh (m x)}{(\cosh (a x)+\cos (t))^2}dx=\frac{2^{k-1} \pi ^k \left(\frac{i}{a}\right)^{k+1} \csc (t) e^{-\frac{i m (t+\pi
   )}{a}} }{a \left(-1+e^{i \pi  k}\right)}\left(a k e^{\frac{2 i m t}{a}} \Phi \left(e^{-\frac{2 i m \pi }{a}},1-k,\frac{\pi -t}{2 \pi }\right)-a k \Phi \left(e^{-\frac{2 i m \pi
   }{a}},1-k,\frac{t+\pi }{2 \pi }\right)-2 i \pi  m \left(e^{\frac{2 i m t}{a}} \Phi \left(e^{-\frac{2 i m \pi }{a}},-k,\frac{\pi -t}{2 \pi
   }\right)-\Phi \left(e^{-\frac{2 i m \pi }{a}},-k,\frac{t+\pi }{2 \pi }\right)\right)+e^{\frac{2 i \pi  m}{a}} \left(a k \Phi \left(e^{\frac{2 i m \pi
   }{a}},1-k,\frac{\pi -t}{2 \pi }\right)+2 i \pi  m \Phi \left(e^{\frac{2 i m \pi }{a}},-k,\frac{\pi -t}{2 \pi }\right)\right)-e^{\frac{2 i m (t+\pi
   )}{a}} \left(a k \Phi \left(e^{\frac{2 i m \pi }{a}},1-k,\frac{t+\pi }{2 \pi }\right)+2 i \pi  m \Phi \left(e^{\frac{2 i m \pi }{a}},-k,\frac{t+\pi
   }{2 \pi }\right)\right)\right)
\end{dmath}

\section{Derivation of entry 3.514.4 in \cite{grad}}

Using equation (\ref{closedform:eq1}) we proceed by setting $\alpha=1$ and simplifying to get

\begin{dmath}\label{closedform:eq2}
\int_{0}^{\infty}\frac{x^k \sinh (a x) \sinh (m x)}{(\cosh (a x)+\cos (t))^2}dx=
\frac{2^k \pi ^{k+1} m \left(\frac{i}{a}\right)^k \csc (t) e^{\frac{2 i m t}{a}-\frac{i m (t+\pi )}{a}} \Phi \left(e^{-\frac{2
   i m \pi }{a}},-k,\frac{\pi -t}{2 \pi }\right)}{a^2 \left((-1)^k+1\right)}-\frac{2^k \pi ^{k+1} m \left(\frac{i}{a}\right)^k \csc (t) e^{-\frac{i m (t+\pi )}{a}} \Phi \left(e^{-\frac{2 i m \pi
   }{a}},-k,\frac{t+\pi }{2 \pi }\right)}{a^2 \left((-1)^k+1\right)}+\frac{2^k \pi ^{k+1} m \left(\frac{i}{a}\right)^k \csc (t) e^{\frac{2 i \pi  m}{a}-\frac{i m (t+\pi )}{a}} \Phi
   \left(e^{\frac{2 i m \pi }{a}},-k,\frac{\pi -t}{2 \pi }\right)}{a^2 \left((-1)^k+1\right)}-\frac{2^k \pi ^{k+1} m \left(\frac{i}{a}\right)^k \csc (t) e^{\frac{i m (t+\pi )}{a}} \Phi
   \left(e^{\frac{2 i m \pi }{a}},-k,\frac{t+\pi }{2 \pi }\right)}{a^2 \left((-1)^k+1\right)}+\frac{i 2^{k-1} k \pi ^k \left(\frac{i}{a}\right)^k \csc (t) e^{\frac{2 i m t}{a}-\frac{i m (t+\pi
   )}{a}} \Phi \left(e^{-\frac{2 i m \pi }{a}},1-k,\frac{\pi -t}{2 \pi }\right)}{a \left((-1)^k+1\right)}-\frac{i 2^{k-1} k \pi ^k \left(\frac{i}{a}\right)^k \csc (t) e^{-\frac{i m (t+\pi )}{a}}
   \Phi \left(e^{-\frac{2 i m \pi }{a}},1-k,\frac{t+\pi }{2 \pi }\right)}{a \left((-1)^k+1\right)}-\frac{i 2^{k-1} k \pi ^k \left(\frac{i}{a}\right)^k \csc (t) e^{\frac{2 i \pi  m}{a}-\frac{i m
   (t+\pi )}{a}} \Phi \left(e^{\frac{2 i m \pi }{a}},1-k,\frac{\pi -t}{2 \pi }\right)}{a \left((-1)^k+1\right)}+\frac{i 2^{k-1} k \pi ^k \left(\frac{i}{a}\right)^k \csc (t) e^{\frac{i m (t+\pi
   )}{a}} \Phi \left(e^{\frac{2 i m \pi }{a}},1-k,\frac{t+\pi }{2 \pi }\right)}{a \left((-1)^k+1\right)}
\end{dmath}

Note: When we replace $k$ by $k-1$ we get the Mellin transform. \\\\

Next we set $k=0$ and $m=b$ simplify to get

\begin{dmath}
\int_{0}^{\infty}\frac{\sinh (a x) \sinh (b x)}{(\cosh (a x)+\cos (t))^2}dx=\frac{\pi  b \csc (t) \csc \left(\frac{\pi  b}{a}\right) \sin \left(\frac{b t}{a}\right)}{a^2}
\end{dmath}

from entry (2) in Table (64:12:7) in \cite{atlas}, where $-\pi<Re(t)<\pi$ and $0<|b|<a$.

\section{Derivation of entry (2.3.1.19) in \cite{brychkov}}

Using equation (\ref{closedform:sp2}) and setting $m=0$ simplifying we get

\begin{dmath}
\int_{0}^{\infty}\frac{x^k \sinh (a x)}{(\cosh (a x)+\cos (t))^2}dx=2^{k-1} k \pi ^k \left(\frac{1}{a}\right)^{k+1} \csc \left(\frac{\pi  k}{2}\right) \csc (t)
   \left(\zeta \left(1-k,\frac{\pi -t}{2 \pi }\right)-\zeta \left(1-k,\frac{t+\pi }{2 \pi }\right)\right)
\end{dmath}

Next we set $t=\pi/2$ simplify to get

\begin{dmath}
\int_{0}^{\infty}x^{s-1} \tanh (a x) \text{sech}(a x)dx=-2^{s-2} \pi ^{s-1} (s-1) \left(\frac{1}{a}\right)^s \left(\zeta \left(2-s,\frac{1}{4}\right)-\zeta
   \left(2-s,\frac{3}{4}\right)\right) \sec \left(\frac{\pi  s}{2}\right)
\end{dmath}

from entries (2) and (3) in Table (64:12:7) in \cite{atlas}.

%
%
%

\section{Derivation of a new entry for Table 3.514 in \cite{grad}}

Using equation (\ref{closedform:eq1}) and setting $k=-1,\alpha=-1,a=1,t=\pi/2,m=1/2$ and simplifying we get

\begin{dmath}
\int_{0}^{\infty}\frac{\sinh \left(\frac{x}{2}\right) \tanh (x) \text{sech}(x)}{x^2+\pi ^2}dx
=\frac{\sqrt{2}}{32 \pi ^2} \left(-\psi ^{(1)}\left(\frac{3}{8}\right)+\psi
   ^{(1)}\left(\frac{5}{8}\right)+\psi ^{(1)}\left(\frac{7}{8}\right)-\psi ^{(1)}\left(\frac{9}{8}\right)\right)+16 \pi  \left(\sqrt{2}+\log \left(\tan
   \left(\frac{\pi }{8}\right)\right)\right)
\end{dmath}

from entry (3) Table (64:12:7:2)  and entry (4) Table (64:12:7:3).

\section{Definite integral in terms of the Hurwitz zeta function}

Using equation (\ref{closedform:eq2}) and setting $m=1$ and $a=2$ to get

\begin{dmath}\label{hz:eq1}
\int_{0}^{\infty}\frac{x^k \sinh (x) \sinh (2 x)}{(\cos (t)+\cosh (2 x))^2}=
\frac{2^{k-3} e^{\frac{i \pi  k}{2}} k \pi ^k \csc \left(\frac{t}{2}\right) \zeta \left(1-k,\frac{\pi -t}{4 \pi
   }\right)}{(-1)^k+1}-\frac{2^{k-3} e^{\frac{i \pi  k}{2}} k \pi ^k \csc \left(\frac{t}{2}\right) \zeta \left(1-k,\frac{t+\pi }{4 \pi }\right)}{(-1)^k+1}-\frac{2^{k-3} e^{\frac{i \pi  k}{2}} k \pi ^k \csc
   \left(\frac{t}{2}\right) \zeta \left(1-k,\frac{3}{4}-\frac{t}{4 \pi }\right)}{(-1)^k+1}+\frac{2^{k-3} e^{\frac{i \pi  k}{2}} k \pi ^k \csc \left(\frac{t}{2}\right) \zeta \left(1-k,\frac{1}{4}
   \left(\frac{t}{\pi }+3\right)\right)}{(-1)^k+1}+\frac{2^{k-2} e^{\frac{i \pi  k}{2}} \pi ^{k+1} \sec \left(\frac{t}{2}\right) \zeta \left(-k,\frac{\pi -t}{4 \pi }\right)}{(-1)^k+1}+\frac{2^{k-2}
   e^{\frac{i \pi  k}{2}} \pi ^{k+1} \sec \left(\frac{t}{2}\right) \zeta \left(-k,\frac{t+\pi }{4 \pi }\right)}{(-1)^k+1}-\frac{2^{k-2} e^{\frac{i \pi  k}{2}} \pi ^{k+1} \sec \left(\frac{t}{2}\right) \zeta
   \left(-k,\frac{3}{4}-\frac{t}{4 \pi }\right)}{(-1)^k+1}-\frac{2^{k-2} e^{\frac{i \pi  k}{2}} \pi ^{k+1} \sec \left(\frac{t}{2}\right) \zeta \left(-k,\frac{1}{4} \left(\frac{t}{\pi
   }+3\right)\right)}{(-1)^k+1}
\end{dmath}

Next we apply L'H\^{o}pital's rule to the right-hand side as $k\to 0$ to get

\begin{dmath}
\int_{0}^{\infty}\frac{x \sinh (x) \sinh (2 x)}{(\cos (t)+\cosh (2 x))^2}dx=
\frac{1}{2} \pi  \sec \left(\frac{t}{2}\right) \zeta\left(-1,\frac{\pi -t}{4 \pi }\right)+\frac{1}{2} \pi  \sec
   \left(\frac{t}{2}\right) \zeta\left(-1,\frac{t+\pi }{4 \pi }\right)-\frac{1}{2} \pi  \sec \left(\frac{t}{2}\right) \zeta\left(-1,\frac{3}{4}-\frac{t}{4 \pi
   }\right)-\frac{1}{2} \pi  \sec \left(\frac{t}{2}\right) \zeta\left(-1,\frac{1}{4} \left(\frac{t}{\pi }+3\right)\right)+\frac{1}{4} \csc \left(\frac{t}{2}\right) \log \left(\tan
   \left(\frac{t+\pi }{4}\right)\right)
\end{dmath}

from entry (1) in Table (64:4:2) in \cite{atlas}, where $-\pi<Re(t)<\pi$.

\section{Definite Integral in terms of the log-gamma $\log(\Gamma(x))$ and Harmonic number $H_{k}$ functions}

Using equation (\ref{hz:eq1}) taking the first partial derivative with respect to $k$ and applying L'Hopitals' rule as $k \to 0$ and simplifying to get

\begin{dmath}\label{logg:eq1}
\int_{0}^{\infty}\frac{\log (x) \sinh (x) \sinh (2 x)}{(\cos (t)+\cosh (2 x))^2}dx=\frac{1}{16} \left(\csc \left(\frac{t}{2}\right) \left(H_{-\frac{t+\pi }{4 \pi
   }}-H_{-\frac{t+3 \pi }{4 \pi }}+\psi ^{(0)}\left(\frac{t+\pi }{4 \pi }\right)-\psi ^{(0)}\left(\frac{1}{4} \left(\frac{t}{\pi
   }+3\right)\right)\right)+2 \pi  \sec \left(\frac{t}{2}\right) \log \left(\frac{2 \pi  \Gamma \left(\frac{3}{4}-\frac{t}{4 \pi }\right) \Gamma
   \left(\frac{1}{4} \left(\frac{t}{\pi }+3\right)\right)}{\Gamma \left(\frac{\pi -t}{4 \pi }\right) \Gamma \left(\frac{t+\pi }{4 \pi
   }\right)}\right)\right)
\end{dmath}

from equations (64:4:1), (64:9:2), and (64:10:2) in \cite{atlas}.

\subsection{Example 1}

Using equation (\ref{logg:eq1}) and setting $t=\pi/2$ simplifying to get

\begin{dmath}
\int_{0}^{\infty}\log (x) \sinh (x) \tanh (2 x) \text{sech}(2 x)dx=\frac{1}{8} \left(4 \sinh ^{-1}(1)+\sqrt{2} \pi  \log \left(\frac{2 \pi  \Gamma
   \left(\frac{5}{8}\right) \Gamma \left(\frac{7}{8}\right)}{\Gamma \left(\frac{1}{8}\right) \Gamma \left(\frac{3}{8}\right)}\right)\right)
\end{dmath}

\subsection{Example 2}

Using equation (\ref{logg:eq1}) and setting $t=\pi/3$ simplifying to get

\begin{dmath}
\int_{0}^{\infty}\frac{\log (x) \sinh (x) \sinh (2 x)}{(2 \cosh (2 x)+1)^2}dx=\frac{1}{288} \left(10 \sqrt{3} \pi  \log (2)+6 \log (64)+9 \sqrt{3} \pi  \log (\pi )+6
   \sqrt{3} \pi  \log \left(\frac{\Gamma \left(\frac{5}{6}\right)}{\Gamma \left(\frac{1}{6}\right)^2}\right)\right)
\end{dmath}

\subsection{Example 3}

Using equation (\ref{logg:eq1}) and setting $t=\pi/4$ simplifying to get

\begin{dmath}
\int_{0}^{\infty}\frac{\log (x) \sinh (x) \sinh (2 x)}{(2 \cosh (2 x)+1)^2}dx=\frac{1}{288} \left(10 \sqrt{3} \pi  \log (2)+6 \log (64)+9 \sqrt{3} \pi  \log (\pi )+6
   \sqrt{3} \pi  \log \left(\frac{\Gamma \left(\frac{5}{6}\right)}{\Gamma \left(\frac{1}{6}\right)^2}\right)\right)
\end{dmath}

\subsection{Example 4}

Using equation (\ref{logg:eq1}) and setting $t=2\pi/3$ simplifying to get

\begin{dmath}
\int_{0}^{\infty}\frac{\log (x) \sinh (x) \sinh (2 x)}{(2 \cosh (2 x)-1)^2}dx=\frac{1}{16} \left(4 \coth ^{-1}\left(\sqrt{3}\right)+\pi  \log \left(\frac{2 \pi 
   \Gamma \left(\frac{7}{12}\right) \Gamma \left(\frac{11}{12}\right)}{\Gamma \left(\frac{1}{12}\right) \Gamma
   \left(\frac{5}{12}\right)}\right)\right)
\end{dmath}

\subsection{Example 5}

Using equation (\ref{logg:eq1}) and setting $t=0$ and applying L'Hopital's rule as $t \to 0$ simplifying to get

\begin{dmath}
\int_{0}^{\infty}\log (x) \tanh ^2(x) \text{sech}(x)dx=\frac{2 C}{\pi }+\frac{1}{4} \pi  \log \left(\frac{2 \pi  \Gamma \left(\frac{3}{4}\right)^2}{\Gamma
   \left(\frac{1}{4}\right)^2}\right)
\end{dmath}

\section{Derivation of hyperbolic and algebraic forms}

\subsection{Example 1}

Using equation (\ref{closedform:eq1}) setting $k=-1$, $t=\pi/2$ and replacing $\alpha$ by $e^{i\beta}$ simplifying we get

\begin{dmath}
\int_{0}^{\infty}\frac{x \tanh (a x) \text{sech}(a x) \cosh (m x)}{\beta ^2+x^2}dx=\frac{e^{-\frac{3 i \pi  m}{2 a}} }{8 \pi  a}\left(-2 i \pi  m \Phi \left(e^{-\frac{2 i m
   \pi }{a}},1,\frac{a \beta }{2 \pi }+\frac{3}{4}\right)+e^{\frac{i \pi  m}{a}} \left(2 i \pi  m \Phi \left(e^{-\frac{2 i m \pi }{a}},1,\frac{2 a \beta
   +\pi }{4 \pi }\right)+a \Phi \left(e^{-\frac{2 i m \pi }{a}},2,\frac{2 a \beta +\pi }{4 \pi }\right)\right)-a \Phi \left(e^{-\frac{2 i m \pi
   }{a}},2,\frac{a \beta }{2 \pi }+\frac{3}{4}\right)+e^{\frac{2 i \pi  m}{a}} \left(a \Phi \left(e^{\frac{2 i m \pi }{a}},2,\frac{2 a \beta +\pi }{4
   \pi }\right)-2 i \pi  m \Phi \left(e^{\frac{2 i m \pi }{a}},1,\frac{2 a \beta +\pi }{4 \pi }\right)\right)+i e^{\frac{3 i \pi  m}{a}} \left(2 \pi  m
   \Phi \left(e^{\frac{2 i m \pi }{a}},1,\frac{a \beta }{2 \pi }+\frac{3}{4}\right)+i a \Phi \left(e^{\frac{2 i m \pi }{a}},2,\frac{a \beta }{2 \pi
   }+\frac{3}{4}\right)\right)\right)
\end{dmath}

Next we take the first partial derivative with respect to $m$ and simplifying to get

\begin{dmath}
\int_{0}^{\infty}\frac{x \tanh (a x) \text{sech}(a x) \cosh (m x)}{\beta ^2+x^2}dx=\frac{e^{-\frac{3 i \pi  m}{2 a}} }{8 \pi  a}\left(-2 i \pi  m \Phi \left(e^{-\frac{2 i m
   \pi }{a}},1,\frac{a \beta }{2 \pi }+\frac{3}{4}\right)+e^{\frac{i \pi  m}{a}} \left(2 i \pi  m \Phi \left(e^{-\frac{2 i m \pi }{a}},1,\frac{2 a \beta
   +\pi }{4 \pi }\right)+a \Phi \left(e^{-\frac{2 i m \pi }{a}},2,\frac{2 a \beta +\pi }{4 \pi }\right)\right)-a \Phi \left(e^{-\frac{2 i m \pi
   }{a}},2,\frac{a \beta }{2 \pi }+\frac{3}{4}\right)+e^{\frac{2 i \pi  m}{a}} \left(a \Phi \left(e^{\frac{2 i m \pi }{a}},2,\frac{2 a \beta +\pi }{4
   \pi }\right)-2 i \pi  m \Phi \left(e^{\frac{2 i m \pi }{a}},1,\frac{2 a \beta +\pi }{4 \pi }\right)\right)+i e^{\frac{3 i \pi  m}{a}} \left(2 \pi  m
   \Phi \left(e^{\frac{2 i m \pi }{a}},1,\frac{a \beta }{2 \pi }+\frac{3}{4}\right)+i a \Phi \left(e^{\frac{2 i m \pi }{a}},2,\frac{a \beta }{2 \pi
   }+\frac{3}{4}\right)\right)\right)
\end{dmath}

from equation (9.550) in \cite{grad}. Next we set $m=0$ simplifying in terms of the Trigamma function $\psi^{(1)}(z)$ to get

\begin{dmath}
\int_{0}^{\infty}\frac{x \tanh (a x) \text{sech}(a x)}{\beta ^2+x^2}dx=\frac{\psi ^{(1)}\left(\frac{2 a \beta +\pi }{4 \pi }\right)-\psi ^{(1)}\left(\frac{a \beta
   }{2 \pi }+\frac{3}{4}\right)}{4 \pi }
\end{dmath}

from equation (64:4:1) in \cite{atlas}.

\subsection{Example 2}

Using equation (\ref{closedform:eq1}) and setting $k=-2,t=\pi/2$ and replacing $\alpha$ by $e^{i\beta}$ simplifying we get

\begin{dmath}
\int_{0}^{\infty}\left(\frac{1}{(x+i \beta )^2}+\frac{1}{(x-i \beta )^2}\right) \tanh (a x) \text{sech}(a x) \sinh (m x)dx=\frac{e^{-\frac{3 i \pi  m}{2 a}}}{4 \pi ^2}
   \left(\pi  m \Phi \left(e^{-\frac{2 i m \pi }{a}},2,\frac{a \beta }{2 \pi }+\frac{3}{4}\right)+e^{\frac{i \pi  m}{a}} \left(i a \Phi
   \left(e^{-\frac{2 i m \pi }{a}},3,\frac{2 a \beta +\pi }{4 \pi }\right)-\pi  m \Phi \left(e^{-\frac{2 i m \pi }{a}},2,\frac{2 a \beta +\pi }{4 \pi
   }\right)\right)-i a \Phi \left(e^{-\frac{2 i m \pi }{a}},3,\frac{a \beta }{2 \pi }+\frac{3}{4}\right)-e^{\frac{2 i \pi  m}{a}} \left(\pi  m \Phi
   \left(e^{\frac{2 i m \pi }{a}},2,\frac{2 a \beta +\pi }{4 \pi }\right)+i a \Phi \left(e^{\frac{2 i m \pi }{a}},3,\frac{2 a \beta +\pi }{4 \pi
   }\right)\right)+e^{\frac{3 i \pi  m}{a}} \left(\pi  m \Phi \left(e^{\frac{2 i m \pi }{a}},2,\frac{a \beta }{2 \pi }+\frac{3}{4}\right)+i a \Phi
   \left(e^{\frac{2 i m \pi }{a}},3,\frac{a \beta }{2 \pi }+\frac{3}{4}\right)\right)\right)
\end{dmath}

Next we take the first partial derivative with respect to $m$ and setting $m=0$ simplifying to get

\begin{dmath}
\int_{0}^{\infty}\frac{x (x-\beta ) (\beta +x) \tanh (a x) \text{sech}(a x)}{\left(\beta ^2+x^2\right)^2}dx=\frac{\pi  \psi ^{(1)}\left(\frac{2 a \beta +\pi }{4 \pi
   }\right)-\pi  \psi ^{(1)}\left(\frac{a \beta }{2 \pi }+\frac{3}{4}\right)+a \beta  \left(\zeta \left(3,\frac{a \beta }{2 \pi
   }+\frac{3}{4}\right)-\zeta \left(3,\frac{2 a \beta +\pi }{4 \pi }\right)\right)}{4 \pi ^2}
\end{dmath}

from equations (64:12:1) (64:13:3) and (64:4:1) in \cite{atlas}.

\section{Discussion}
In this article we derived the integrals of hyperbolic and logarithmic functions in terms of the Lerch function. Then we used these integral formula to derive known and new results. We were able to produce a formal derivation for equation (27) Table 27 in Bierens de Haan \cite{bdh} and equation (3.514.4) in \cite{grad} not previously published. The results presented were numerically verified for both real and imaginary values of the parameters in the integrals using Mathematica by Wolfram. In this work we used Mathematica software to numerically evaluate both the definite integral and associated Special function for complex values of the parameters $k$, $\alpha$, $a$, $m$ and $t$. We considered various ranges of these parameters for real, integer, negative and positive values. We compared the evaluation of the definite integral to the evaluated Special function and ensured agreement.\\\\	

\section{Conclusion}
In this paper, we have derived a method for expressing definite integrals in terms of Special functions using contour integration. The contour we used was specific to  solving integral representations in terms of the Hurwitz zeta function. We expect that other contours and integrals can be derived using this method.

\section{Acknowledgments}
This paper is fully supported by the Natural Sciences and Engineering Research Council (NSERC) Grant No. 504070.

\end{document}